\documentclass[a4paper,12pt]{article}

\usepackage{amssymb}
\usepackage{tipa}

\usepackage[T2A]{fontenc}
\usepackage[cp1251]{inputenc}
\usepackage{amsmath}

\usepackage[dvips]{graphicx}
\usepackage[
russian, 
english
]{babel}

\usepackage{wasysym}

\def \GK{\operatorname{GK}}
\def\Pid{\operatorname{Pid}}
\def\Var{\operatorname{Var}}
\def\Ht{\operatorname{Ht}}
\def\PI{{\rm PI}}

\newtheorem{theorem}{Теорема}[section]
\newtheorem{corollary}[theorem]{Следствие}
\newtheorem{lemma}[theorem]{Лемма}
\newtheorem{proposition}[theorem]{Предложение}
\newtheorem{notation}[theorem]{Обозначение}
\newtheorem{remark}[theorem]{Замечание}

\newtheorem{algorithm}[theorem]{Алгоритм}
\newtheorem{ques}[theorem]{Вопрос}

\newtheorem{definitionhead}[theorem]{Определение}
\newenvironment{definition}{\begin{definitionhead}%
\sl}{\end{definitionhead}}

\begin{document}


\title{
\ \hbox to \textwidth{\normalsize      
16R+05E\hfill} 
\ \hbox to \textwidth{\normalsize
УДК 512.5+512.64+519.1\hfill}\\[1ex]
Subexponential estimations in Shirshov's Height theorem}

\author{
Alexei Belov, Mikhail Kharitonov}

\maketitle

\begin{abstract}
In 1993 E.~I.~Zelmanov asked the following question in Dniester Notebook:

``Suppose that $F_{2, m}$ is a 2-generated associative ring with the identity $x^m=0.$ Is it true, that the nilpotency degree of $F_{2, m}$ has exponential growth?$"$

We show that the nilpotency degree of $l$-generated associative algebra with the identity $x^d=0$ is smaller than $\Psi(d,d,l),$ where $$\Psi(n,d,l)=2^{18} l (nd)^{3 \log_3 (nd)+13}d^2.$$
We give the definitive answer to E.~I.~Zelmanov by this result. It is the consequence of one fact, which is based on combinatorics of words. Let $l$, $n$ и $d \geqslant n$ be positive integers. Then all the words over
alphabet of cardinality  $l$ which length is greater than
 $\Psi(n,d,l)$ are either $n$-divided or contain $d$-th power of subword, where a word  $W$ is {\it $n$-divided}, if it can be represented in the
following form  $W=W_0W_1\cdots W_n$ such that $W_1\prec
W_2\prec\dots\prec W_n$. The symbol $\prec$ means lexicographical order here. A.~I.~Shirshov
proved that the set of non $n$-divided words over alphabet of
cardinality $l$ has bounded height $h$ over the set $Y$ consisting
of all the words of degree $\leqslant n-1$. Original Shirshov's estimation was just recursive, in 1982 double
exponent was obtained by A.~G.~Kolotov and in 1993 A.~Ya.~Belov obtained exponential
estimation. We show, that $h<\Phi(n,l)$, where
$$\Phi(n,l) = 2^{87} l\cdot n^{12\log_3 n + 48}.$$
Our proof uses Latyshev idea of Dilworth theorem application.
\bigskip
\bigskip

В 1993 году Е.~И.~Зельманов поставил следующий вопрос в Днестровской тетради:

``Пусть $F_{2,m}$ -- свободное $2$-порождённое ассоциативное кольцо с тождеством $x^m=0.$ Верно ли, что класс нильпотентности кольца $F_{2,m}$ растёт экспоненциально по $m?"$

В данной работе показано, что в $l$-порождённой ассоциативной алгебре с тождеством $x^d=0$ класс нильпотентности меньше, чем $\Psi(d,d,l),$ где $$\Psi(n,d,l)=2^{18} l (nd)^{3 \log_3 (nd)+13}d^2.$$
Тем самым получен окончательный ответ на вопрос Е.~И.~Зельманова. Данный результат является следствием следующего факта, относящегося к комбинаторике слов. Пусть $l$, $n$ и $d\geqslant n$ -- некоторые натуральные числа. Тогда все
слова над $l$-буквенным алфавитом длины не меньше, чем
$\Psi(n,d,l)$, либо содержат $x^d$, либо являются
$n$-разбиваемыми, где слово $W$
называется {\it $n$-разбиваемым}, если его можно представить в
виде $W=W_0W_1\cdots W_n,$ где подслова $W_1,\dots,W_n$ идут в
порядке лексикографического убывания.

Из не $n$-разбиваемых  слов
состоит базис алгебры с тождеством степени $n$. А.~И.~Ширшов
показал, что множество слов, не являющихся $n$-разбиваемыми, над
алфавитом из $l$ букв имеет ограниченную высоту $h$ над  $Y$~--
множеством слов степени не выше $n-1$.
Изначальная оценка числа $h$ у Ширшова носила рекурсивный характер, в 1982 А.~Г.~Колотов
году получил двойную экспоненту, а в 1993 году А.~Я.~Белов получил экспоненциальную оценку. Мы показываем, что
$h<\Phi(n,l)$, где
$$\Phi(n,l) = 2^{87} l\cdot n^{12\log_3 n + 48}.$$
Доказательство использует идею В.~Н.~Латышева, связанную с
применением теоремы Дилуорса к исследованию не
$n$-раз\-би\-ва\-ем\-ых слов.
\end{abstract}



\medskip


{\bf Ключевые слова:}Теорема Ширшова о высоте,
комбинаторика слов, $n$-разбиваемый, теоремы Дилуорса, проблемы Бернсайдовского типа.

\section{Введение}
\subsection{Теорема Ширшова о высоте}

В 1958 году А.~И.~Ширшов доказал свою знаменитую теорему о высоте (\cite{Shirshov1}, \cite{Shirshov2}).

\begin{definition}
Назовем слово $W$ {\em $n$-разбиваемым}, если $W$ можно
представить в виде $W = vu_1u_2\cdots u_n$ так чтобы $u_1\succ
u_2\succ\dots\succ u_n$.
\end{definition}

В этом случае при любой нетождественной перестановке $\sigma$
подслов $u_i$ получается слово $W_\sigma =
vu_{\sigma(1)}u_{\sigma(2)}\cdots u_{\sigma(n)}$,
лексикографически меньшее $W$. Это свойство некоторые авторы берут
за основу определения понятия $n$-разбиваемости.

\begin{definition}
Назовём \PI-алгебру $A$ алгеброй {\em ограниченной высоты
$\Ht_Y(A)$ над множеством слов $Y = \{ u_1, u_2,\ldots\}$}, если
$h$ -- минимальное число такое, что любое слово $x$ из $A$ можно
представить в виде
$$x = \sum_i \alpha_i u_{j_{(i,1)}}^{k_{(i,1)}}
u_{j_{(i, 2)}}^{k_{(i,2)}}\cdots
u_{j_{(i,r_i)}}^{k_{(i,r_i)}},
$$
причем $\{r_i\}$ ограничены числом $h$ в совокупности. Множество
$Y$ называется {\em базисом Ширшова} для $A$.
\end{definition}

\medskip
 {\bf Теорема Ширшова о высоте.} (\cite{Shirshov1}, \cite{Shirshov2})
 {\it Множество всех не $n$-разбиваемых слов в конечно порождённой алгебре с допустимым полиномиальным тождеством имеет ограниченную
 высоту $H$ над множеством слов степени не выше $n-1$.
 }
 \medskip

Из теоремы о высоте вытекает решение ряда проблем теории колец
(см. ниже). Проблемы Бернсайдовского типа, связанные с теоремой о высоте, рассмотрены в обзоре \cite{Zelmanov}. Однако, авторы убеждены, что теорема Ширшова о высоте
является фундаментальным фактом комбинаторики слов безотносительно приложений к $\PI$-теории. (Все наши доказательства
элементарны и не выходят за рамки работы со словами.) К сожалению,
специалисты по комбинаторике данный факт должным образом пока не
оценили. Что касается самого понятия {\it $n$-разбиваемости}, то
оно также представляется фундаментальным. Оценки, полученные
В.~Н.~Латышевым на $\xi_n(k)$ -- количество не $n$-разбиваемых
полилинейных слов от $k$ символов, привели к фундаментальным
результатам в $\PI$-теории. Вместе с тем, это количество суть не
что иное, как количество расстановок чисел от $1$ до $k$ таких,
что никакие $n$ из них (не обязательно стоящие подряд) не идут в
порядке убывания. То же самое -- перечисление всех
перестановочно упорядоченных множеств диаметра $n$. (Множество
называется {\it перестановочно упорядоченным}, если его порядок
есть пересечение двух линейных порядков.)

 Из этой теоремы вытекает положительное решение проблем Бернсайдовского типа
 для $\PI$-алгебр.
В самом деле, пусть в ассоциативной алгебре над полем выполняется
полиномиальное тождество $f(x_1,\ldots,x_n)=0$. Можно доказать,
что тогда в ней выполняется и допустимое полилинейное тождество
(т.е. полиномиальное тождество, у которого хотя бы один
коэффициент при членах высшей степени равен единице):
$$
x_1 x_2\cdots x_n = \sum_{\sigma}\alpha_{\sigma}x_{\sigma(1)}
x_{\sigma(2)}\cdots x_{\sigma(n)},$$
 где $\alpha_{\sigma}$ принадлежат основному полю. В этом случае, если
 $W = v u_1 u_2 \cdots u_n$ является $n$-разбиваемым, то для любой
 перестановки $\sigma$ слово $W_{\sigma} = vu_{\sigma(1)}u_{\sigma(2)}\cdots u_{\sigma(n)}$
 лексикографически меньше слова $W$, т.е. $n$-разбиваемое слово
 можно представить в виде линейной комбинации лексикографически меньших слов.
 Значит, $\PI$-алгебра имеет базис из не $n$-разбиваемых слов.
 В силу теоремы Ширшова о высоте,
 $\PI$-алгебра имеет ограниченную высоту.
 Как следствие имеем, что если в $\PI$-алгебре выполняется тождество $x^n = 0$, то эта
 алгебра -- нильпотентна, т.е. все ее слова длины больше, чем некоторое $N$,
 тождественно равны 0. Обзоры, посвященные теореме о высоте, содержатся в работах
 \cite{BBL97,Kem09,BelovRowenShirshov}.

Из этой теоремы вытекает положительное решение проблемы Куроша, и
других проблем бернсайдовского типа для $\PI$-колец. Ведь если
$Y$~-- базис Ширшова, и все элементы из $Y$~-- алгебраичны, то
алгебра $A$ конечномерна. Тем самым теорема Ширшова дает явное
указание множества элементов, алгебраичность которых ведет к
конечномерности всей алгебры. Из этой теоремы вытекает

\begin{corollary}[Berele]
Пусть $A$~-- конечно порожденная $\PI$-алгебра. Тогда
$$\GK(A)<\infty.$$
\end{corollary}

$\GK(A)$ -- это {\it размерность Гельфанда--Кириллова алгебры $A$},
т.е.
$$\GK(A)=\lim_{n\to\infty}\ln V_A(n)/\ln(n),$$
где $V_A(n)$ есть {\it функция роста алгебры $A$}, т.е.
размерность векторного пространства, порожденного словами степени
не выше $n$ от образующих $A$.

В самом деле, достаточно заметить, что число решений неравенства\\
$k_1 |v_1|+\cdots+k_h|v_h|\leqslant n$, где $h\leqslant H$, превосходит $N^{H}$, и потому $\GK(A)\leqslant \Ht(A)$.

Число $m=\deg(A)$ будет обозначать {\it степень алгебры} или
минимальную степень тождества, которое в ней выполняется.
$n=\Pid(A)$ есть {\it сложность} алгебры $A$ или максимальное $k$
такое, что ${\Bbb M}_k$ --- алгебра матриц размера $k$, принадлежит
многообразию $\Var(A)$, порожденному алгеброй $A$.

Вместо понятия {\it высоты} удобнее пользоваться близким понятием
{\it существенной высоты}.

\begin{definition}
Алгебра $A$ имеет {\em существенную высоту $h=H_{Ess}(A)$} над
конечным множеством $Y$, называемым {\em $s$-базисом алгебры $A$}, если можно
выбрать такое конечное множество $D\subset A$, что $A$ линейно
представима элементами вида $t_1\cdot\ldots\cdot t_l$, где $l\leqslant
2h+1$, и $\forall i (t_i\!\in\! D \vee t_i=y_i^{k_i};y_i\in Y)$,
причем множество таких $i$, что $t_i\not\in D$, содержит не более
$h$ элементов. Аналогично определяется {\em существенная высота}
множества слов.
\end{definition}

Говоря неформально, любое длинное слово есть произведение
периодических частей и ``прокладок'' ограниченной длины.
Существенная высота есть число таких периодических кусков, а
обычная еще учитывает ``прокладки''.

В связи с теоремой о высоте возникли следующие вопросы:

\begin{enumerate}

\item На какие классы колец можно распространить теорему о высоте?

\item Над какими $Y$ алгебра $A$ имеет ограниченную высоту? В
частности, какие наборы слов можно взять в качестве $\{v_i\}$?

\item Как устроен вектор степеней $(k_1,\ldots,k_h)$? Прежде
всего: какие множества компонент этого вектора являются
существенными, т.е. какие наборы $k_i$ могут быть одновременно
неограниченными? Какова существенная высота? Верно ли что
множество векторов степеней обладает теми или иными свойствами
регулярности?

\item Как оценить высоту?
\end{enumerate}

Перейдем к обсуждению поставленных вопросов.
\subsection{Неассоциативные обобщения}
 Теорема о высоте была
распространена на некоторые классы колец, близких к ассоциативным.
С.~В.~Пчелинцев \cite{Pchelintcev} доказал ее для альтернативного
и $(-1,1)$ случаев, С.~П.~Мищенко \cite{Mishenko1} получил аналог
теоремы о высоте для алгебр Ли с разреженным тождеством. В работе
автора \cite{Belov1} теорема о высоте была доказана для некоторого
класса колец, асимптотически близких к ассоциативным, куда
входят, в частности, альтернативные и йордановы $\PI$-алгебры.
\subsection{Базисы Ширшова}
Пусть $A$ -- $\PI$-алгебра, и подмножество
$M\subseteq A$ является ее $s$-базисом. Тогда, если все
элементы множества $M$ алгебраичны над $K$, то алгебра $A$
конечномерна (проблема Куроша).  Ограниченность существенной
высоты над $Y$ влечет ``положительное решение проблемы Куроша над
$Y$''. Обратное утверждение менее тривиально.

\begin{theorem}  [А.~Я.~Белов]         \label{ThKurHmg}
а) Пусть $A$~-- градуированная $\PI$-алгебра, $Y$~-- конечное
множество однородных элементов. Тогда если при всех $n$ алгебра
$A/Y^{(n)}$ нильпотентна, то $Y$ есть $s$-базис $A$. Если при этом
$Y$ порождает $A$ как алгебру, то $Y$~-- базис Ширшова алгебры
$A$.

б) Пусть $A$~--- $\PI$-алгебра, $M\subseteq A$~--- некоторое курошево
подмножество в $A$. Тогда $M$~--- $s$-базис алгебры $A$.
\end{theorem}

$Y^{(n)}$ обозначает идеал, порожденный $n$-ыми степенями элементов
из $Y$. Множество $M\subset A$ называется {\it курошевым}, если
любая проекция $\pi\colon A\otimes K[X]\to A'$, в которой образ
$\pi(M)$ цел над $\pi(K[X])$, конечномерна над $\pi(K[X])$.
Мотивировкой этого понятия служит следующий пример.  Пусть
$A={\Bbb Q}[x,1/x]$. Любая проекция $\pi$ такая, что $\pi(x)$
алгебраичен, имеет конечномерный образ. Однако множество $\{x\}$
не является $s$-ба\-зис\-ом алгебры ${\Bbb Q}[x,1/x]$. Таким
образом, ограниченность существенной высоты есть некоммутативное
обобщение свойства {\it целости}.
\subsection{Базисы Ширшова, состоящие из слов}
Описание базисов
Ширшова, состоящих из слов, дает следующая теорема:

\begin{theorem}[\cite{BBL97}, \cite{BR05}]            \label{ThBelheight}
Множество слов $Y$ является базисом Ширшова алгебры $A$ тогда и
только тогда, когда для любого слова $u$ длины не выше $m =
\Pid(A)$~--- сложности алгебры $A$~--- множество $Y$ содержит
слово, циклически сопряженное к некоторой степени слова $u$.
\end{theorem}
\subsection{Существенная высота}
Ясно, что размерность
Гель\-фан\-да--Ки\-рил\-ло\-ва оценивается существенной высотой и
что $s$-базис является базисом Ширшова, тогда и только тогда,
когда он порождает $A$ как алгебру. В представимом случае имеет
место и обратное утверждение.

\begin{theorem}[А.~Я.~Белов, \cite{BBL97}]
Пусть $A$~-- конечно порожденная представимая алгебра и пусть
$H_{Ess}{}_Y(A)<\infty$. Тогда $H_{Ess}{}_Y(A)=\GK(A)$.
\end{theorem}

\begin{corollary}[В.~Т.~Марков]
Размерность Гель\-фан\-да--Ки\-рил\-ло\-ва ко\-неч\-но
по\-рож\-ден\-ной представимой алгебры есть целое число.
\end{corollary}

\begin{corollary}
Если $H_{Ess}{}_Y(A)<\infty$, и алгебра $A$ представима, то
$H_{Ess}{}_Y(A)$ не зависит от выбора $s$-базиса $Y$.
\end{corollary}

В этом случае размерность Гель\-фан\-да---Ки\-рил\-ло\-ва также
равна существенной высоте в силу локальной представимости относительно свободных алгебр.

\medskip
{\bf Строение векторов степеней.} Хотя в представимом случае
размерность Гель\-фан\-да-Ки\-рил\-ло\-ва и существенная высота
ведут себя хорошо. Тем не менее даже тогда множество векторов
степеней может быть устроено плохо~--- а именно, может быть
дополнением к множеству решений системы
экс\-по\-нен\-ци\-аль\-но-по\-ли\-но\-ми\-аль\-ных диофантовых
уравнений \cite{BBL97}. Вот почему существует пример представимой
алгебры с трансцендентным рядом Гильберта. Однако для относительно
свободной алгебры ряд Гильберта рационален \cite{Belov501}.
\medskip
\subsection{$n$-разбиваемость и теорема Дилуорса}
Значение понятия {\it
$n$-раз\-би\-ва\-ем\-ос\-ти} выходит за рамки проблематики,
относящейся к проблемам бернсайдовского типа. Оно играет роль и
при изучении полилинейных слов, в оценке их количества, где {\it
полилинейным} называется слово, в которое каждая буква входит
не более одного раз. В.~Н.~Латышев применил теорему Дилуорса для
получения оценки числа не являющихся $m$-разбиваемыми полилинейных слов степени $n$, над
алфавитом $\{a_1,\dots,a_n\}$. Эта
оценка:  ${(m - 1)}^{2n}$ и она близка к реальности. Напомним эту
теорему.

\medskip
{\bf Теорема Дилуорса}: {\it Пусть $n$ -- наибольшее количество
элементов антицепи данного конечного частично упорядоченного
множества $M$. Тогда $M$ можно разбить на $n$  попарно
непересекающихся цепей.}
\medskip

Рассмотрим полилинейное слово $W$ из $n$ букв. Положим $a_i\succ
a_j$, если $i>j$ и буква $a_i$ стоит в слове $W$ правее $a_j$.
Условие не $m$-разбиваемости означает отсутствие антицепи из $m$
элементов. Тогда по теореме Дилуорса все позиции (и,
соответственно, буквы $a_i$) разбиваются на $(m-1)$ цепь. Сопоставим
каждой цепи свой цвет. Тогда раскраска позиций и раскраска букв
однозначно определяет слово $W$. А число таких раскрасок не
превосходит $(m-1)^n\times (m-1)^n=(m-1)^{2n}$.

Из данной оценки следует выполнимость полилинейных тождеств,
отвечающих неприводимому модулю, диаграмма Юнга которого содержит
квадрат $n^4$. Это в свою очередь, во-первых, позволило получить
прозрачное доказательство теоремы Регева о том, что тензорное
произведение $\PI$-алгебр снова является $\PI$-алгеброй,
во-вторых, установить существование разреженного тождества в общем
случае, а также тождества Капелли в конечно порожденном случае (тем
самым, в частности, доказать теорему о нильпотентности радикала),
и в-третьих, осуществить ``супертрюк'' А.~Р.~Кемера, сводящий
изучение тождеств общих алгебр к изучению супер-тождеств конечнопорождённых супералгебр в нулевой характеристике. Смежные вопросы рассмотрены в работе \cite{BP07,Lot83,02}.

Вопросы, связанные с перечислением полилинейных слов, не
являющихся $n$-раз\-би\-ва\-е\-мы\-ми, имеют самостоятельный
интерес. (Например, существует биекция между не
$3$-раз\-би\-ва\-е\-мы\-ми словами и числами Каталана.) С одной
стороны, это чисто комбинаторная задача, с другой стороны, она
связана с рядом коразмерностей для алгебры общих матриц.
Исследование полилинейных слов представляется чрезвычайно важным.
В.~Н.~Латышев  (см. например \cite{LatyshevMulty}) поставил
проблему конечной базируемости множества старших полилинейных слов
для $T$-идеала относительно взятия надслов и изотонных
подстановок. Из этой проблемы вытекает проблема Шпехта для
полилинейных многочленов, имеется тесная связь с проблемой слабой
нётеровости групповой алгебры бесконечной финитарной
симметрической группы над полем положительной характеристики (для
нулевой характеристики это было установлено А.~Залесским). Для
решения проблемы Латышева надо уметь переводить свойства
$T$-идеалов на язык полилинейных слов. В работах
\cite{BBL97,Belov1} была попытка осуществить программу перевода
структурных свойств алгебр на язык комбинаторики слов. На язык
полилинейных слов такой перевод осуществить проще, в дальнейшем
можно получить информацию и о словах общего вида.

В данной работе мы переносим  технику В.~Н.~Латышева на
не полилинейный случай, что позволяет получить субэкспоненциальную
оценку в теореме Ширшова о высоте. Г.~Р.~Челноков предложил идею этого переноса в 1996 году.
\subsection{Оценки высоты}
Первоначальное доказательство А.~И.~Ширшова
хотя и было чисто комбинаторным (оно основывалось на технике
элиминации, развитой им в алгебрах Ли, в частности, в
доказательстве теоремы о свободе), однако оно давало только
примитивно рекурсивные оценки. Позднее А.~Т.~Колотов
\cite{Kolotov} получил оценку на $\Ht(A)\leqslant l^{l^n}$\
($n=\deg(A)$,\, $l$~-- число образующих). А.~Я.~Белов в работе \cite{Bel92} показал, что $\Ht(n,l)<2nl^{n+1}$.  Экспоненциальная оценка теоремы Ширшова о высоте изложена также в работах \cite{BR05}, \cite{Dr00}, \cite{Kh11}, \cite{Kh11(3)}. Данные оценки улучшались в работах А.~Клейна \cite{Klein,Klein1}. В 2001 году Е.~С.~Чибриков в работе \cite{Ch01} доказал, что \\ $\Ht(4,l) \geqslant  (7k^2-2k).$ М.~И.~Харитонов получил в работах \ref{Kh11}, \cite{Kh11(2)} верхние и нижние оценки на существенную высоту. В 2011 году А.~А.~Лопатин \cite{Lop11} получил следующий результат:

\begin{theorem}
Пусть $C_{n,l}$ -- степень нильпотентности свободной $l$-порождённой алгебры и удовлетворяющей тождеству $x^n=0.$ Пусть $p$ -- характеристика базового поля алгебры -- больше чем $n/2.$ Тогда
$$(1): C_{n,l}<4\cdot 2^{n/2} l.$$
\end{theorem}
По определению $C_{n,l}\leqslant \Psi(n, n, l).$
Заметим, что для малых $n$ оценка  (1) меньше, чем полученная в данной работе оценка $\Psi(n, n, l),$ но при росте $n$ оценка $\Psi(n,n,l)$ асимптотически лучше оценки (1).

Е.~И.~Зельманов поставил следующий вопрос в Днестровской тетради в 1993 году:
\begin{ques}
Пусть $F_{2,m}$ -- свободное $2$-порождённое ассоциативное кольцо с тождеством $x^m=0.$ Верно ли, что класс нильпотентности кольца $F_{2,m}$ растёт экспоненциально по $m?$
\end{ques}

Наша работа отвечает на вопрос Е.~И.~Зельманова следующим образом: в действительности искомый класс нильпотентности растёт субэкспоненциально.
\subsection{Полученные результаты}
{\bf Основной результат} работы состоит в следующем:

\begin{theorem}     \label{c:main2}
Высота множества не $n$-разбиваемых слов над $l$-буквенным
алфавитом относительно множества слов длины меньше $n$ не
превышает $\Phi(n,l)$, где
$$\Phi(n,l) = E_1 l\cdot n^{E_2+12\log_3 n} ,$$ где
$E_1 = 4^{21\log_3 4 + 17}, E_2 = 30\log_3 4 + 10.$
\end{theorem}
Из данной теоремы путем некоторого огрубления и упрощения оценки получается, что
при фиксированном $l$ и $n \rightarrow\infty$
$$\Phi(n,l) < 2^{87} l\cdot n^{12\log_3 n + 48}
= n^{12(1+o(1))\log_3{n}},$$

а при фиксированном $n$ и $l\rightarrow\infty$
$$\Phi(n,l) < C(n)l.$$

\begin{corollary}
Высота $l$-порождённой $\PI$-алгебры с допустимым полиномиальным
тождеством степени $n$ над множеством слов длины меньше $n$ не
превышает $\Phi(n,l)$.
\end{corollary}

Кроме того, доказывается субэкспоненциальная оценка, которая лучше при малых $n$:

\begin{theorem}     \label{t1:log_2}
Высота множества не $n$-разбиваемых слов над $l$-буквенным
алфавитом относительно множества слов длины меньше $n$ не
превышает $\Phi(n,l)$, где
$$\Phi(n,l) = 2^{40} l\cdot n^{38+8\log_2 n}.$$
\end{theorem}

Как следствие получаются субэкспоненциальные оценки на индекс
нильпотентности $l$-порожденных ниль-алгебр степени $n$ для
произвольной характеристики.

Другим основным результатом нашей работы является следующая

\begin{theorem}      \label{c:main1}
Пусть $l$, $n$ и $d\geqslant n$ -- некоторые натуральные числа. Тогда все
$l$-порождённые слова длины не меньше, чем $\Psi(n,d,l)$, либо
содержат $x^d$, либо являются $n$-разбиваемыми, где
$$
\Psi(n,d,l)=4^{5+3\log_3 4} l (nd)^{3 \log_3 (nd)+(5+6\log_3 4)}d^2.
$$
\end{theorem}
Из данной теоремы путем некоторого огрубления и упрощения оценки получается, что
при фиксированном $l$ и $nd \rightarrow\infty$
$$\Psi(n,d,l) < 2^{18} l (nd)^{3 \log_3 (nd)+13}d^2
= (nd)^{3(1+o(1))\log_3(nd)},$$

а при фиксированном $n$ и $l\rightarrow\infty$
$$\Psi(n,d,l) < C(n,d)l.$$

\begin{corollary}
 Пусть $l$, $d$ -- некоторые натуральные числа.
Пусть в ассоциативной  $l$-порождённой алгебре $A$ выполнено тождество
$x^{d}=0$. Тогда ее индекс нильпотентности меньше, чем
$\Psi(d,d,l)$.
\end{corollary}

Кроме того, доказывается субэкспоненциальная оценка, которая лучше при малых $n$ и $d$:

\begin{theorem}     \label{t2:log_2}
Пусть $l$, $n$ и $d\geqslant n$ -- некоторые натуральные числа. Тогда все
$l$-порождённые слова длины не меньше, чем $\Psi(n,d,l)$, либо
содержат $x^d$, либо являются $n$-разбиваемыми, где
$$
\Psi(n,d,l) = 256 l(nd)^{2\log_2 (nd)+10}d^2.$$
\end{theorem}

\begin{notation}
Для вещественного числа $x$ положим $\ulcorner x\urcorner := -[-x].$ Таким образом мы округляем нецелые числа в большую сторону.
\end{notation}
В процессе доказательства теоремы \ref{c:main2} доказывается следующая теорема, оценивающая существенную высоту:

\begin{theorem} \label{ThThick}
Существенная высота $l$-порождённой $PI$-алгебры с допустимым полиномиальным тождеством степени $n$ над множеством слов длины меньше $n$ меньше, чем $\Upsilon (n, l),$ где
$$\Upsilon (n, l) = 2n^{3\ulcorner\log_3 n\urcorner + 4} l.$$
\end{theorem}

В работе \cite{Bog01} установлено, что индекс нильпотентности
$l$-порожденного ниль-полукольца степени $n$ совпадает с индексом
нильпотентности $l$-порожденного ниль-кольца степени $n$, причем
сложение не обязательно коммутативно. (Там же приведены примеры
не нильпотентных ниль почти колец индекса $2$.) Таким образом,
наши результаты распространяются и на случай полуколец.
\medskip

\subsection{О нижних оценках}
Сравним полученные результаты  с нижней
оценкой для высоты. Высота алгебры $A$ не меньше ее размерности
Гель\-фан\-да--Ки\-рил\-ло\-ва $\GK(A)$. Для алгебры
$l$-порождённых общих матриц порядка $n$ данная размерность, равна
$(l-1)n^2+1$ (см. \cite{Procesi} а также \cite{Bel04}). В то же
время, минимальная степень тождества этой алгебры равна $2n$ в
силу теоремы Ам\-и\-цу\-ра--Ле\-виц\-ко\-го. Имеет место следующее

\begin{proposition}
Высота $l$-порожденной $\PI$-алгебры степени $n$, а также
множества  не $n$-раз\-би\-ва\-ем\-ых слов над $l$-бук\-вен\-ным
алфавитом, не менее, чем $(l-1)n^2/4+1$.
\end{proposition}

Нижние оценки на индекс нильпотентности были установлены
Е.~Н.~Кузьминым в работе \cite{Kuz75}.  Е.~Н.~Кузьмин привел
пример $2$-порожденной алгебры с тождеством $x^n=0$, индекс
нильпотентности которой строго больше $(n^2+n-2)/2$. Вопрос нахождения нижних оценок рассматривается в работах \cite{Kh11}, \cite{Kh11(2)}.

В то же время для случая нулевой характеристики и счетного числа
образующих Ю.~П.~Размыслов (см. например, \cite{Razmyslov3})
получил верхнюю оценку на индекс нильпотентности равную $n^2$.

Вначале мы докажем теорему \ref{c:main1}, в следующей главе мы
займемся оценками существенной высоты, т.е. количества различных
периодических фрагментов в не $n$-разбиваемом слове.
\medskip

{\bf Благодарности.} Авторы признательны В.~Н.~Латышеву,
А.~В.~Михалёву и всем участникам семинара ``Теория колец'' за
внимание к работе, а также участникам семинара МФТИ под руководством А.~М.~Райгородского.

\section{Оценки на появление степеней подслов}
\subsection{План доказательства теоремы \ref{c:main1}}

В леммах \ref{Lm0.1}, \ref{c:lem1.2} и \ref{c:lem1.3} описываются достаточные условия для присутствия периода длины $d$ в не $n$-разбиваемом слове $W$. В лемме \ref{c:lem1.4} связываются понятия $n$-разбиваемости слова $W$ и множества его хвостов. После этого определённым образом выбирается подмножество множества хвостов слова $W$, для которого можно применить теорему Дилуорса. Затем мы раскрашиваем хвосты и их первые буквы в соответствии с принадлежностью к цепям, полученным при применении теоремы Дилуорса.

Необходимо изучить, в какой позиции начинают отличаться соседние хвосты в каждой цепи. Вызывает интерес, с какой $``$частотой$"$ эта позиция попадёт в $p$-хвост для некоторого $p\leqslant n$. Потом мы несколько обобщаем наши рассуждения, деля хвосты на сегменты по несколько букв, а затем рассматривая, в какой сегмент попала позиция, в которой начинают отличаться друг от друга соседние хвосты в цепи.
В лемме \ref{c:lem2} связываются рассматриваемые $``$частоты$"$ для $p$-хвостов и $kp$-хвостов для $k = 3$.

В завершение доказательства строится иерархическая структура на основе применения леммы \ref{c:lem2}, т.~е. рассматриваем сначала сегменты $n$-хвостов, потом подсегменты этих сегментов и т.~д. Далее мы рассматриваем наибольшее возможное количество хвостов из подмножества, для которого была применена теорема Дилуорса, после чего оцениваем сверху общее количество хвостов, а, значит, и букв слова $W$.

\subsection{Свойства периодичности и $n$-разбиваемости}

\smallskip


Пусть $a_1, a_2,\ldots ,a_l$ -- алфавит, над которым проводится построение слов. Порядок $a_1\prec a_2\prec\dots\prec a_l$ индуцирует
лексикографический порядок на словах над заданным алфавитом. Для удобства введём следующие определения:

\begin{definition}
а) Если в слове $v$ содержится подслово вида $u^t,$ то будем
говорить, что в слове $v$ содержится период длины $t.$

б) Если слово $u$ является началом слова $v$, то такие слова
называют {\em несравнимыми}.

в) Слово $v$ -- {\em хвост} слова $u$, если найдется слово $w$
такое, что $u=wv$.

г) Слово $v$ -- {\em $k$-хвост} слова $u$, если $v$ состоит из
$k$ первых букв некоторого хвоста $u$.

г*) {\em $k$-начало} то же самое, что и $k$-хвост.

д) Пусть слово $u$ {\em левее} слова $v$, если начало слова $u$ левее начала слова $v$.
\end{definition}

\begin{notation}
а) Для вещественного числа $x$ положим $\ulcorner x\urcorner := -[-x].$

б) Обозначим как $|u|$ длину слова $u$.
\end{notation}

Для доказательства потребуются следующие достаточные условия наличия периода:

\begin{lemma}       \label{Lm0.1}
В слове $W$ длины $x$ либо первые $[x/d]$ хвостов попарно
сравнимы, либо в слове $W$ найдется период длины $d$.
\end{lemma}

$\RHD$  Пусть в слове $W$ не нашлось слова вида $u^{d}$. Рассмотрим первые
$[x/{d}]$ хвостов. Предположим, что среди них нашлись 2
несравнимых хвоста $v_1$ и $v_2$. Пусть $v_1=u\cdot v_2$. Тогда
$v_2=u\cdot v_3$ для некоторого $v_3$. Тогда $v_1=u^2\cdot v_3$.
Применяя такие рассуждения, получим, что $v_1=u^{d}\cdot
v_{{d}+1}$, так как $|u|<x/ {d}$, $|v_2|\geqslant ({d}-1)x/ {d}$.
Противоречие.$\LHD$

\begin{lemma}     \label{c:lem1.2}
Если в слове $V$ длины $k\cdot t$ не больше $k$ различных подслов
длины $k$, то $V$ включает в себя период длины $t$.
\end{lemma}

$\RHD$ Докажем лемму индукцией по $k$. База при $k = 1$ очевидна. Если
найдется не больше, чем $(k - 1)$ различных подслов длины $(k -
1)$, то применяем индукционное предположение. Если существуют $k$
различных подслов длины $(k - 1)$, то каждое подслово длины $k$
однозначно определяется своими первыми $(k - 1)$ буквами. Значит,
$V=v^t$, где $v$ -- $k$-хвост $V$.$\LHD$

\begin{definition}
а) Слово $W$ -- \textit{$n$-разбиваемо в обычном смысле}, если
найдутся $u_1, u_2,\ldots,u_n$ такие, что $W=v\cdot u_1\cdots
u_n$, при этом $u_1\succ\ldots \succ u_n$.

б) В текущем доказательстве слово $W$ будем называть
\textit{$n$-разбиваемым в хвостовом смысле}, если найдутся хвосты
$u_1,\ldots,u_n$ такие, что $u_1\succ u_2\succ \ldots \succ u_n$ и
для любого $i=1, 2,\ldots, n - 1$ начало $u_i$ слева от начала
$u_{i+1}$. Если особо не оговорено противное, то под
\textit{$n$-разбиваемыми} словами мы подразумеваем $n$-разбиваемые
в хвостовом смысле.

в) Слово $W$ -- \textit{$n$-сократимое}, если оно либо
$n$-разбиваемо в обычном смысле, либо найдется слово вида $u^{d}\subseteq W$.
\end{definition}

Теперь опишем достаточное условие $n$-сократимости и его связь с $n$-разбиваемостью.

\begin{lemma}          \label{c:lem1.3}
Если в слове $W$ найдутся $n$ одинаковых непересекающихся подслов
$u$ длины $n\cdot{d}$, то $W$ -- $n$-сократимое.
\end{lemma}

$\RHD$ Предположим противное. Рассмотрим хвосты
$u_1,u_2,\ldots,u_n$ слова $u$, которые начинаются с каждой из его
первых $n$ букв. Перенумеруем хвосты так, чтобы выполнялись
неравенства: $u_1~\succ\ldots\succ~u_n.$ Из леммы 1 они
несравнимые. Рассмотрим подслово $u_1,$ лежащее в самом левом экземпляре слова $u,$ подслово $u_2$ -- во втором слева,$\ldots, u_n$ -- в $n$-ом слева. Получили
$n$-разбиение слова $W$. Противоречие.$\LHD$

\begin{lemma}      \label{c:lem1.4}
Если слово $W$ является $4nd$-разбиваемым, то оно --
$n$-сократимое.
\end{lemma}

$\RHD$ Предположим противное. Рассмотрим порядковые номера позиций букв $a_i$,
где $a_1<a_2<\ldots<a_{4nd}$, с которых начинаются хвосты $u_i$, разбивающие $W$. Положим $a_{4nd+1} = |W|$. Если $W$ -- не $n$-сократимое, то найдётся такое число $i$, что $1\leqslant i\leqslant 4(n-1)d + 1,$ что для любых $i\leqslant b<c\leqslant d<e\leqslant i+4d$ $(a_c - a_b)$-хвост $u_b$ и $(a_e-a_d)$-хвост $u_d$ -- несравнимы. Сравним числа $a_{i+2d} - a_i$ и $a_{i+4d} - a_{i+2d}$. Можно считать, что $a_{i+4d} - a_{i+2d}\geqslant a_{i+2d} - a_i$. Пусть $a_{j+1} - a_j = \inf\limits_k {(a_{k+1} - a_k)}, 0\leqslant j<2d.$ Сравним числа $j$ и $d$. Можно считать, что $j<d.$ По предположению $(a_{2d}-a_j)$-хвост $u_j$ и $(a_{2d} - a_{j+1})$-хвост $u_{j+1}$ несравнимы с $(a_{4d}-a_{2d})$-хвостом $u_{2d}$. Т. к. $a_{4d}-a_{2d}\geqslant a_{2d}-a_j>a_{2d}-a_{j+1},$ то $(a_2d-a_j)$-хвост $u_j$ и $(a_{2d}-a_{j+1})$-хвост $u_{j+1}$ несравнимы между собой. Так как ${{a_{2d}-a_j}\over {a_{2d}-a_{j+1}}}\leqslant {{d+1}\over {d}},$ то $(a_{j+1}-a_j)$-хвост $u_j$ в степени $d$ содержится в $(a_2d-a_j)$-хвосте $u_j$.
Противоречие. $\LHD$

\begin{corollary}
Если слово $W$ -- не $n$-разбиваемо в обычном смысле, то $W$ не
$4nd$-разбиваемо (в хвостовом смысле).
\end{corollary}

\begin{notation}
Положим $p_{n, d}:=4nd-1$.
\end{notation}


Пусть $W$ -- не $n$-сократимое слово. Рассмотрим $U$ -- $[\left|
W\right|/d]$-хвост слова $W$. Тогда $W$ -- не $(p_{n, d}+1)$-разбиваемое.
Пусть $\Omega$ -- множество хвостов слова $W$, которые начинаются
в $U$. Тогда лемме \ref{Lm0.1} любые два элемента из $\Omega$
сравнимы. Естественным образом строится биекция между $\Omega$,
буквами $U$ и натуральными числами от $1$ до
$\left|\Omega\right|=\left|U\right|$.

Введем слово $\theta$ такое, что $\theta$ лексикографически меньше
любого слова.

\begin{remark}
В текущем доказательстве теоремы \ref{c:main1} все хвосты мы
предполагаем лежащими в $\Omega$.
\end{remark}

\section{Оценки на появление периодических фрагментов}

\paragraph{Применение теоремы Дилуорса.} Для хвостов $u$ и $v$ положим
$u<v$, если $u \prec v$ и, кроме того, $u$ левее $v$. Тогда по
теореме Дилуорса $\Omega$ можно разбить на $p_{n, d}$ цепей, где в
каждой цепи $u \prec v$, если $u$ левее $v$. Покрасим начальные
позиции хвостов в $p_{n, d}$ цветов в соответствии с принадлежностью к
цепям. Фиксируем натуральное число $p$. Каждому
натуральному числу $i$ от 1 до $\left|\Omega\right|$ сопоставим
$B^p(i)$ -- упорядоченный набор из $p_{n, d}$ слов $\{f(i,j)\}$
построенных по следующему правилу:

{\it Для каждого $j = 1, 2,\ldots, p_{n, d}$\ положим

$$f(i,j)=\left\{\max \
f\leqslant i: f\ \mbox{раскрашено\ в\ цвет}~j\right\}.$$

Если такого $f$ не найдется, то слово из $B^p(i)$ на позиции
$j$ считаем равным $\theta$, в противном случае это слово считаем
равным $p$-хвос\-ту, который начинается с $f(i,j)$-ой
буквы. }

Неформально говоря, мы наблюдаем, с какой скоростью хвосты "эволюционируют'' в своих цепях, если рассматривать последовательность позиций слова $W$ как ось времени.

\subsection{Наборы $B^p(i)$, процесс на позициях}

\begin{lemma}[О процессе]   \label{c:lem}
Дана последовательность $S$ длины $|S|$, составленная из слов
длины $(k-1)$. Каждое из них состоит из $(k-2)$ символа $``0"$ и одной
$``1"$. Пусть $S$ удовлетворяет следующему условию:

{\em если для некоторого $0 < s \leqslant k-1$ найдутся $p_{n, d}$
слов, в которых $``1"$ стоит на $s$-ом месте, то между первым и
$p_{n, d}$-ым из этих слов найдется слово, в котором $``1"$ стоит
строго меньше, чем на $s$-ом месте}; $L(k-1)=\sup\limits_S |S|$.

Тогда $L(k-1)\leqslant p_{n, d}^{k-1}-1$.
\end{lemma}

$\RHD$ $L(1) \leqslant p_{n, d}-1$. Пусть $L(k-1) \leqslant {p_{n, d}^{k-1}} -
1$. Покажем, что $L(k) \leqslant {p_{n, d}^{k}} - 1$.
Рассмотрим слова, у которых символ $``1"$ стоит на первом месте.
Их не больше $p_{n, d}-1$. Между любыми двумя из них, а также перед
первым и после последнего, количество слов не больше $L(k-1)
\leqslant {p_{n, d}^{k-1}} - 1$. Получаем, что
$$
L(k) \leqslant p_{n, d} - 1 +
\left(p_{n, d}\right)\left({\left(p_{n, d}\right)^{k-1}}-1\right) =
{\left(p_{n, d}\right)^{k}}-1\LHD$$

Нам требуется ввести некоторую величину, которая бы численно оценивала скорость "эволюции'' наборов $B^p(i)$:

\begin{definition}
Положим $$\psi(p):= \left\{\max \ k:
B^p(i)=B^p(i+k-1)\right\}.$$
В частности по лемме \ref{c:lem1.2},
$\psi(p_{n,d})\leqslant p_{n, d} d$.

Для заданного $\alpha$ определим разбиение последовательности первых $\left|\Omega\right|$ позиций
${i}$ слова $W$ на классы эквивалентности $\AC_\alpha$ следующим образом: $i
\AC_\alpha j$, если $B^\alpha(i)=B^\alpha(j).$
\end{definition}

\begin{proposition}
Для любых натуральных $a<b$ имеем $\psi(a) \leqslant \psi(b).$
\end{proposition}

\begin{lemma}[Основная] \label{c:lem2}
Для любых натуральных чисел $a, k$ верно неравенство
$$\psi(a)\leqslant p_{n,d}^k\psi(k\cdot a)+k\cdot a$$
\end{lemma}

$\RHD$
Рассмотрим по наименьшему представителю из каждого класса $\AC_{k\cdot a}$. Получена последовательность позиций $\{i_j\}.$ Теперь рассмотрим все $i_j$ и $B^{k\cdot a}(i_j)$ из одного класса эквивалентности по $\AC_a.$ Пусть он состоит из $B^{k\cdot a}(i_j)$ при $i_j\in[b, c).$ Обозначим за $\{i_j\}'$ отрезок последовательности $\{i_j\},$ для которого $i_j\in[b, c-k\cdot a).$

Фиксируем некоторое натуральное число $r, 1\leqslant r\leqslant p_{n,d}.$ Назовём все $k\cdot a$-начала цвета $r$, начинающиеся с позиций слова $W$ из $\{i_j\}'$, представителями типа $r$. Все представители типа $r$, будут попарно различны, так как они начинаются с наименьших позиций в классах эквивалентности по $\AC_{k\cdot a}.$ Разобьём каждый представитель типа $r$ на $k$ сегментов длины $a$. Пронумеруем сегменты внутри каждого представителя типа $r$ слева направо числами от нуля до $(k-1)$. Если найдутся $(p_{n,d}+1)$ представителей типа $r$, у которых совпадают первые $(t-1)$ сегментов, но которые попарно различны в $t$-ом, где $t$ -- натуральное число, $1\leqslant t\leqslant k-1$, то найдутся две первых буквы $t$-го сегмента одного цвета. Тогда позиции, с которых начинаются эти сегменты, входят в разные классы эквивалентности по $\AC_a.$

Применим лемму \ref{c:lem} следующим образом: во всех представителях типа $r$, кроме самого правого, будем считать сегменты {\it единичными}, если именно в них находится наименьшая позиция, в которой текущий представитель типа $r$ отличается от предыдущего. Остальные сегменты считаем {\it нулевыми.}

Теперь можно применить лемму о процессе с параметрами, совпадающими с заданными в условии леммы. Получаем, что в последовательности $\{i_j\}'$ будет не более $p_{n,d}^{k-1}$ представителей типа $r$. Тогда в последовательности $\{i_j\}'$ , будет не более $p_{n,d}^{k}$ членов. Таким образом, $c-b\leqslant p_{n,d}^k\psi(k\cdot a)+k\cdot a.$ $\LHD$

\subsection{Завершение доказательств теорем \ref{c:main1} и \ref{t2:log_2}}

Пусть
$$a_0 = 3^{\ulcorner \log_3 p_{n,d}\urcorner}, a_1 = 3^{\ulcorner \log_3 p_{n,d}\urcorner-1},\ldots,a_{{\ulcorner \log_3 p_{n,d}\urcorner}} =1.$$
При этом $\left|W\right|\leqslant d\left|\Omega\right| + d$ в силу леммы \ref{Lm0.1}.

Так как набор $B^1(i)$ принимает не более $(1+p_{n,d}l)$ различных значений, то $\left|W\right|\leqslant d(1+p_{n,d}l)\psi(1) + d.$

По лемме \ref{c:lem2}
$$\psi(1)< (p_{n,d}^3 + p_{n,d})\psi(3)<(p_{n,d}^3 + p_{n,d})^2\psi(9)<\cdots <(p_{n,d}^3 + p_{n,d})^{\ulcorner \log_3 p_{n,d}\urcorner}\psi(p_{n,d})\leqslant
$$
$$\leqslant (p_{n,d}^3 + p_{n,d})^{\ulcorner \log_3 p_{n,d}\urcorner}p_{n,d}d
$$

Подставляя $p_{n,d} = 4nd-1,$ получаем
$$\left|W\right|< 4^{5+3\log_3 4} l (nd)^{3 \log_3 (nd)+(5+6\log_3 4)}d^2.
$$

Отсюда имеем утверждение {\bf теоремы \ref{c:main1}.}

Доказательство теоремы \ref{t2:log_2} завершается также, только вместо последовательности
$$a_0 = 3^{\ulcorner \log_3 p_{n,d}\urcorner}, a_1 = 3^{\ulcorner \log_3 p_{n,d}\urcorner-1},\ldots,a_{{\ulcorner \log_3 p_{n,d}\urcorner}} =1$$
рассматривается последовательность
$$a_0 = 2^{\ulcorner \log_2 p_{n,d}\urcorner}, a_1 = 2^{\ulcorner \log_2 p_{n,d}\urcorner-1},\ldots,a_{{\ulcorner \log_2 p_{n,d}\urcorner}} =1.$$

\section{Оценка существенной высоты.}
В данном разделе мы продолжаем доказывать основную теорему \ref{c:main2}.
Попутно доказывается теорема \ref{ThThick}. Будем смотреть на позиции букв слова $W$ как на ось времени, то есть подслово $u$ встретилось раньше подслова $v$, если $u$ целиком лежит левее $v$ внутри слова $W$.



\subsection{Вычленение различных периодических фрагментов в слове~$W$} \label{c:sub1}

Обозначим за $s$ количество подслов слова $W$ с периодом длины меньше $n$, в которых период повторяется больше $2n$ раз и которые попарно разделены сравнимыми с предыдущим периодом подсловами длины больше $n$. Пронумеруем их от начала к концу слова: $x_1^{2n}, x_2^{2n},\ldots,x_s^{2n}$. Таким образом $W=y_0x_1^{2n}y_1x_2^{2n}\cdots x_s^{2n}y_s.$



Если найдётся $i$ такое, что слово $x_i$ длины не меньше $n$, то в слове $x_i^2$ найдутся $n$ попарно сравнимых хвостов, а, значит, слово $x_i^{2n}$-- $n$-разбиваемое. Получаем, что число $s$ не меньше, чем существенная высота слова $W$ над множеством слов длины меньше $n.$

\begin{definition}
Слово $u$ назовем {\em нециклическим}, если $u$ нельзя представить
в виде $v^k$, где $k>1$.
\end{definition}

\begin{definition}
{\em Слово-цикл $u$} -- слово $u$ со всеми его сдвигами по циклу.
\end{definition}

\begin{definition}
Слово $W$ называется {\em сильно $n$-разбиваемым}, если его
можно представить в виде $W=W_0W_1\cdots W_n$, где подслова
$W_1,\dots,W_n$ идут в порядке лексикографического убывания, и
каждое из слов $W_i, i=1, 2,\ldots, n$ начинается с некоторого
слова $z_i^k\in Z$, все $z_i$ различны.
\end{definition}

\begin{lemma}\label{lem4.10}
Если найдётся число $m, 1\leqslant m<n,$ такое, что существуют $(2n-1)$ попарно несравнимых слов длины $m: x_{i_1},\ldots ,x_{i_{2n-1}},$ то $W$ -- $n$-разбиваемое.
\end{lemma}
$\RHD$ Положим $x:=x_{i_1}.$ Тогда в слове $W$ найдутся непересекающиеся подслова\\ $x^{p_1}v'_1,\ldots ,x^{p_{2n-1}}v'_{2n-1},$ где $p_1,\ldots ,p_{2n-1}$ -- некоторые натуральные числа, большие $n,$ а $v'_1,\ldots ,v'_{2n-1}$ -- некоторые слова длины $m,$ сравнимые с $x, v'_1=v_{i_1}.$ Тогда среди слов $v'_1,\ldots ,v'_{2n-1}$ найдутся либо $n$ лексикографически больших $x$, либо $n$ лексикографически меньших $x$. Можно считать, что $v'_1,\ldots ,v'_n$ -- лексикографически больше $x$. Тогда в слове $W$ найдутся подслова $v'_1, xv'_2,\ldots ,x^{n-1}v'_n,$ идущие слева направо в порядке лексикографического убывания.$\LHD$

Рассмотрим некоторое число $m, 1\leqslant m<n.$ Разобьём все $x_i$ длины $m$ на эквивалентности по сильной несравнимости и выберем по одному представителю из каждого класса эквивалентности. Пусть это слова $x_{i_1},\ldots ,x_{i_s'},$
где $s'$ -- некоторое натуральное число. Так как подслова $x_i$ являются периодами, будем рассматривать их как слово-циклы.

\begin{notation}

$v_k := x_{i_k}$

Пусть $v(k, i)$, где $i$ -- натуральное число от 1 до $m$, -- циклический сдвиг слова $v_k$ на $(k - 1)$ позиций вправо, то есть $v(k, 1) = v_k$, а первая буква слова $v(k, 2)$ является второй буквой слова $v_k$. Таким образом, $\{ v(k, i)\}_{i=1}^m$ -- слово-цикл слова $v_k$. Заметим, что для любых $1\leqslant i_1, i_2\leqslant p, 1\leqslant j_1, j_2\leqslant m$ слово $v(i_1, j_1)$ сильно несравнимо со словом $v(i_2, j_2)$.
\end{notation}


\begin{remark}
Случаи $m = 2, 3, n-1$ рассмотрены в работах \cite{Kh11},\cite{Kh11(2)},\cite{Kh11(3)}.
\end{remark}

\subsection{Применение теоремы Дилуорса}      \label{c:sub2}

Рассмотрим множество $\Omega' = \{ v(i, j)\}$, где $1\leqslant i\leqslant p, 1\leqslant j\leqslant m.$ Введём следующий порядок на словах $v(i, j):$

$v(i_1, j_1)\succ v(i_2, j_2),$ если

1) $v(i_1, j_1)> v(i_2, j_2)$

2) $i_1 > i_2$

\begin{lemma}\label{c:lem4}
Если в множестве $\Omega'$ для порядка $\succ$ найдётся антицепь длины $n$, то слово $W$ будет $n$-разбиваемым.
\end{lemma}
$\RHD$ Пусть нашлась антицепь длины $n$ из слов $v(i_1, j_1), v(i_2, j_2),\ldots,v(i_n, j_n);$\\ $i_1\leqslant i_2\leqslant\cdots\leqslant i_n.$ Если все неравенства между $i_k$ -- строгие, то слово $W$ -- $n$-разбиваемое по определению.

Предположим, что для некоторого числа $r$ нашлись $i_{r+1} = \cdots = i_{r+k}$, где либо $r = 0$, либо $i_r < i_{r+1}$. Кроме того, $k$ -- такое натуральное число, что
либо $k = n - r$, либо $i_{r+k} < i_{r+k+1}$.

Слово $s_{i_{r+1}}$ -- периодическое, следовательно, оно представляется в виде произведения $n$ экземпляров слова $v^2_{i_{r+1}}$. $v^2_{i_{r+1}}$ содержит слово-цикл $v_{i_{r+1}}$. Значит, в слове $s_{i_{r+1}}$ можно выбрать непересекающиеся подслова, идущие в порядке лексикографического убывания, равные $v(i_{r+1}, j_{r+1}),\ldots,v(i_{r+k}, j_{r+k})$ соответственно. Таким же образом поступаем со всеми множествами равных индексов  в последовательности $\{i_r\}_{r=1}^n$. Получаем $n$-разбиваемость слова $W$. Противоречие.$\LHD$

Значит, множество $\Omega'$ можно разбить на $(n-1)$ цепь.

\begin{notation}
Положим $q_n = (n-1)$.
\end{notation}

\subsection{Наборы $C^\alpha(i)$, процесс на позициях}

Покрасим первые буквы слов из $\Omega'$ в $q_n$ цветов в соответствии с
принадлежностью к цепям. Покрасим также числа от $1$ до $\left|\Omega '
\right|$ в соответствующие цвета. Фиксируем натуральное число
$\alpha\leqslant m$. Каждому числу $i$ от 1 до $\left|\Omega '
\right|$ сопоставим упорядоченный набор слов $C^\alpha(i)$,
состоящий из $q_n$ слов по следующему правилу:

\medskip
{\it Для каждого $j=1,2,\ldots,q_n$\ положим

$f(i,j)=\{\max \ f\leqslant i:$ существует $k$ такое, что $v(f,k)$ раскрашено в цвет $j$ и
$\alpha$-хвост, который начинается с $f$, состоит только из
букв, являющихся первыми буквами хвостов из $\Omega '\}$.

Если такого $f$ не найдется, то слово из $C^\alpha(i)$ считаем
равным $\theta$, в противном случае это слово считаем равным
$\alpha$-хвосту, слова $v(f,k)$. }

\begin{notation}
Положим  \\$\phi(a)=\{\max\ k:$ для некоторого $ i$ верно $C^a(i)=C^a(i+k-1)\}$.

Для заданного $a\leqslant m$ определим разбиение последовательности
слово-циклов $\{i\}$ слова $W$ на классы эквивалентности следующим
образом: $i \AC_a j,$ если $C^a(i) = C^a(j)$.
\end{notation}

Заметим, что построенная конструкция во многом аналогична построенной в доказательстве теоремы \ref{c:main1}. Обращаем внимание на схожесть $B^a(i)$ и $C^a(i)$, а также $\psi(a)$ и $\phi(a)$.

\begin{lemma}\label{lem:m}
$\phi(m) \leqslant q_n/m$.
\end{lemma}

$\RHD$ Напомним, что слово-циклы были пронумерованы. Рассмотрим слово-циклы с номерами $i, i + 1,\ldots ,i + [q_n/m].$ Ранее было показано, что каждый слово-цикл состоит из $m$ различных слов. Рассмотрим теперь слова в слово-циклах $i, i + 1,\ldots ,i + [q_n/m]$ как элементы множества $\Omega'.$ При таком рассмотрении у первых букв из слово-циклов появляются свои позиции. Всего рассматриваемых позиций не меньше $n.$ Следовательно, среди них найдутся две позиции одного цвета. Тогда в силу сильной несравнимости слово-циклов имеем утверждение леммы.$\LHD$

\begin{proposition}
Для любых натуральных $a<b$ имеем $\phi(a) \leqslant \phi(b).$
\end{proposition}

\begin{lemma}[Основная] \label{lem:thick}
Для натуральных чисел $a, k$, таких, что $ak\leqslant m$ верно неравенство
$$\phi(a)\leqslant p_{n,d}^k\phi(k\cdot a)$$
\end{lemma}

$\RHD$
Рассмотрим по наименьшему представителю из каждого класса $\AC_{k\cdot a}$. Получена последовательность позиций $\{i_j\}.$ Теперь рассмотрим все $i_j$ и $C^{k\cdot a}(i_j)$ из одного класса эквивалентности по $\AC_a.$ Пусть он состоит из $C^{k\cdot a}(i_j)$ при $i_j\in[b, c).$ Обозначим за $\{i_j\}'$ отрезок последовательности $\{i_j\},$ для которого $i_j\in[b, c).$

Фиксируем некоторое натуральное число $r, 1\leqslant r\leqslant q_n.$ Назовём все $k\cdot a$-начала цвета $r$, начинающиеся с позиций слова $W$ из $\{i_j\}'$, представителями типа $r$. Все представители типа $r$, будут попарно различны, так как они начинаются с наименьших позиций в классах эквивалентности по $\AC_{k\cdot a}.$ Разобьём каждый представитель типа $r$ на $k$ сегментов длины $a$. Пронумеруем сегменты внутри каждого представителя типа $r$ слева направо числами от нуля до $(k-1)$. Если найдутся $(q_n+1)$ представителей типа $r$, у которых совпадают первые $(t-1)$ сегментов, но которые попарно различны в $t$-ом, где $t$ -- натуральное число, $1\leqslant t\leqslant k-1$, то найдутся две первых буквы $t$-го сегмента одного цвета. Тогда позиции, с которых начинаются эти сегменты, входят в разные классы эквивалентности по $\AC_a.$

Применим лемму \ref{c:lem} следующим образом: во всех представителях типа $r$, кроме самого правого, будем считать сегменты {\it единичными}, если именно в них находится наименьшая позиция, в которой текущий представитель типа $r$ отличается от предыдущего. Остальные сегменты считаем {\it нулевыми.}

Теперь мы можем применить лемму о процессе с параметрами, совпадающими с заданными в условии леммы. Получаем, что в последовательности $\{i_j\}'$ будет не более $q_n^{k-1}$ представителей типа $r$. Тогда в последовательности $\{i_j\}'$ , будет не более $q_n^{k}$ членов. Таким образом, $c-b\leqslant q_n^k\phi(k\cdot a).$ $\LHD$

\subsection{Завершение доказательства теоремы \ref{ThThick}}\label{ThThick:end}

Пусть
$$a_0 = 3^{\ulcorner \log_3 p_{n,d}\urcorner}, a_1 = 3^{\ulcorner \log_3 p_{n,d}\urcorner-1},\ldots,a_{{\ulcorner \log_3 p_{n,d}\urcorner}} =1.$$
Подставляя эти $a_i$ в леммы
\ref{lem:thick} и \ref{lem:m}, получаем, что
$$\phi(1)\leqslant q_n^3\phi(3)\leqslant q_n^9\phi(9)\leqslant\cdots \leqslant q_n^{3\ulcorner \log_3 m\urcorner}\phi(m)\leqslant
$$
$$\leqslant q_n^{3\ulcorner \log_3 m\urcorner+1}. $$

Так как $C_i^{1}$ принимает не более $1+q_n l$ различных значений, то

$$\left|\Omega'\right|< q_n^{3\ulcorner \log_3 m\urcorner+1}(1+q_n l)< n^{3\ulcorner \log_3 n\urcorner+2} l.$$

По лемме \ref{lem4.10} получаем, что количество $x_i$ длины $m$ меньше $2n^{3\ulcorner \log_3 n\urcorner+3} l.$

Имеем, что количество всех $x_i$ меньше $2n^{3\ulcorner \log_3 n\urcorner+4} l.$

То есть $s<2n^{3\ulcorner \log_3 n\urcorner+4} l.$ Таким образом, теорема \ref{ThThick} доказана.

\section{Доказательство теоремы \ref{t1:log_2} и основной теоремы~\ref{c:main2}}

\subsection{План доказательства}
Будем снова под {\em $n$-разбиваемым
словом} подразумевать {\em $n$-раз\-би\-ва\-е\-мое} в обычном
смысле.
Сначала мы находим необходимое количество фрагментов с длиной периода не меньше $2n$ в слове $W$. Это можно сделать, просто разбив слово $W$ на подслова большой длины, к которым применяется теорема \ref{c:main1}. Однако мы можем улучшить оценку, если сначала выделим в слове $W$ периодический фрагмент с длиной периода не менее $4n$, затем рассмотрим $W_1$ -- слово $W$ с ``вырезанным'' периодическим фрагментом $u_1$. У слова $W_1$ выделяем фрагмент с длиной периода не менее $4n$, после чего рассматриваем $W_2$ -- слово $W_1$ с ``вырезанным'' периодическим фрагментом $u_2$. У слова $W_2$ так же вырезаем периодический фрагмент. Далее продолжаем этот процесс, подробнее описанный в алгоритме \ref{c:al}. Затем мы по вырезанным фрагментам восстанавливаем первоначальное слово $W$. После этого показывается, что в слове $W$ подслово $u_i$ чаще всего не является произведением большого количества не склеенных подслов. В лемме \ref{c:lem3} доказывается, что применение алгоритма \ref{c:al} дает необходимое количество подслов слова $W$ с длиной периода не меньше $2n$ среди вырезанных подслов.

\subsection{Суммирование существенной высоты и степени нильпотентности}

\begin{notation}
$\Ht(w)$ -- высота слова $w$ над множеством слов степени не выше $n$.
\end{notation}
Рассмотрим слово $W$ с высотой $\Ht(W)>\Phi(n,l)$. Теперь для него
проведём следующий алгоритм:

\begin{algorithm}  \label{c:al}
{\ }\\

\begin{description}
    \item[Первый шаг.] По теореме \ref{c:main1} в слове $W$ найдётся подслово с длиной периода
    $4n$. Пусть $W_0=W=u'_1x_{1'}^{4n}y'_1$, причём слово $x_{1'}$ -- нециклическое.
    Представим $y'_1$ в виде $y'_1=x_{1'}^{r_2}y_1$, где $r_2$ -- максимально возможное
    число. Слово $u'_1$ представим как $u'_1=u_1x_{1'}^{r_1}$, где $r_1$ -- наибольшее возможное. Обозначим за $f_1$ следующее слово:
$$W_0=u_1x_{1'}^{4n+r_1+r_2}y_1=u_1f_1y_1.$$
    Назовём позиции, входящие в слово $f_1$, {\em скучными,} последнюю позицию слова $u_1$ -- {\em скучной типа $1$,} вторую с конца позицию $u_1$ -- {\em скучной типа $2$,} и так далее, $n$-ую с конца позицию $u_1$ -- {\em скучной типа $n$.} Положим $W_1 = u_1 y_1.$

    \item[$k$-ый шаг.] Рассмотрим слова $u_{k-1},\ y_{k-1},\ W_{k-1}=u_{k-1}y_{k-1}$,
    построенные на предыдущем шаге. Если $|W_{k-1}|\geqslant\Phi(n,l),$ то применим теорему \ref{c:main1} к слову $W$ с тем условием, что процесс в основной лемме \ref{c:lem2} будет вестись только по не скучным позициям и скучным позициям типа больше $ka,$ где $k$ и $a$ -- параметры леммы \ref{c:lem2}.

    Таким образом, в слове $W_{k-1}$ найдётся
нециклическое подслово с длиной периода $4n$, так что
$$W_{k-1}=u'_kx_{k'}^{4n}y'_k.$$
При этом положим
$$r_1 := \sup \{r: u'_k = u_k x_{k'}^r\},\quad r_2 :=\sup
\{r: y'_k = x_{k'}^r y_k\}.$$

(Отметим, что слова в наших рассуждениях могут быть пустыми.)\\
Определим  $f_k$ из равенства:
$$W_{k-1}= u_kx_{k'}^{4n+r_1+r_2}y_k = u_kf_ky_k.$$
Назовём позиции, входящие в слово $f_k$, {\em скучными}, последнюю позицию слова $u_k$ -- {\em скучной типа $1$,} вторую с конца позицию $u_k$ -- {\em скучной типа $2$,} и так далее, $n$-ую с конца позицию $u_k$ -- {\em скучной типа $n$}. Если позиция в процессе алгоритма определяется как скучная двух типов, то будем считать её скучной того типа, который меньше. Положим $W_k = u_k y_k.$
\end{description}

\end{algorithm}

\begin{notation}
Проведём $4t+1$ шагов алгоритма \ref{c:al}. Рассмотрим
первоначальное слово $W$. Для каждого натурального $i$ из отрезка
$[1,4t]$ имеет место равенство
$$
W = w_0f_i^{(1)}w_1f_i^{(2)}\cdots f_i^{(n_i)} w_{n_i}
$$
для некоторых подслов $w_j$. Здесь $f_i = f_i^{(1)}\cdots
f_i^{(n_i)}$. Также мы считаем, что при $1\leqslant j\leqslant
n_i-1$ подслово $w_j$ -- непустое. Пусть $s(k)$ -- количество
индексов $i\in[1,4t]$ таких, что $n_i = k$.
\end{notation}

Для доказательства теоремы \ref{c:main1} требуется найти как можно больше длинных периодических фрагментов. Помочь в этом сможет следующая лемма:

\begin{lemma}   \label{c:lem3}
$s = s(1) + s(2) \geqslant 2t$.
\end{lemma}

$\RHD$ Назовём {\it монолитным}  подслово $U$ слова $W$, если
\begin{enumerate}
    \item $U$ является произведением слов вида $f_i^{(j)}$,
    \item $U$ не является подсловом слова, для которого выполняется предыдущее
свойство (1).
\end{enumerate}

Пусть после $(i-1)$-го шага алгоритма \ref{c:al} в слове $W$
содержится $k_{i-1}$ монолитных подслов. Заметим, что $k_i
\leqslant k_{i-1}-n_i+2$.

 Тогда если $n_i \geqslant 3$, то $k_i\leqslant k_{i-1}-1.$
Если же $n_i\leqslant 2$, то $k_i\leqslant k_{i-1}+1$. При этом
$k_1=1$, $k_t \geqslant 1 = k_1$. Лемма доказана. $\LHD$

\begin{corollary} \label{c:cor}
$$\sum\limits_{k=1}^\infty {k\cdot s(k)} \leqslant 10t\leqslant 5s.(\ref{c:cor})$$
\end{corollary}

$\RHD$ Из доказательства леммы \ref{c:lem3} получаем, что
$\sum\limits_{n_i\geqslant 3} {(n_i-2)} \leqslant 2t$.

По определению
$\sum\limits_{k=1}^\infty {s(k)} =4t$, т.е.
$\sum\limits_{k=1}^\infty {2s(k)} =8t$.

Складывая эти два неравенства и применяя лемму \ref{c:lem3}, получаем доказываемое
неравенство \ref{c:cor}.$\LHD$

\begin{proposition}
Высота слова $W$ будет не больше $$\Psi(n,4n,l)+\sum\limits_{k=1}^\infty {k\cdot s(k)}\leqslant \Psi(n,4n,l)+ 5s.$$
\end{proposition}

Далее будем рассматривать только $f_i$ с $n_i \leqslant 2$.

\begin{notation}
Если $n_i = 1$, то положим $f'_i := f_i$.

Ежели $n_i = 2$, то положим $f'_i := f_i^{(j)}$, где $f_i^{(j)}$
-- слово с наибольшей длиной между $f_i^{(1)}$ и $f_i^{(2)}$.

Слова $f'_i$ упорядочим в соответствии с их близостью к началу
$W$. Получим последовательность $f'_{m_1},\ldots ,f'_{m_s}$, где
$s'=s(1)+s(2)$, положим $f''_i := f'_{m_i}$. Пусть $f''_i = w'_i
x_{i''}^{p_{i''}}w''_i$, где хотя бы одно из слов $w'_i, w''_i$ --
пустое.
\end{notation}

\begin{remark}  \label{c:pr}
Можно считать, что мы первыми шагами алгоритма \ref{c:al}
выбрали все те $f_i$, для которых $n_i =1$.
\end{remark}

Теперь рассмотрим $z'_j$ -- подслова $W$
следующего вида:
$$z'_j = x_{(2j-1)''}^{p_{(2j-1)''}+\gimel}v_j,
\gimel\geqslant 0, |v_j| = |x_{(2j-1)''}|,$$

при этом $v_j$ не равно $x_{(2j - 1)''}$, начало $z'_j$ совпадает
с началом периодического подслова в $f''_{2j-1}$. Покажем, что
$z'_j$ не пересекаются.

В самом деле. Если $f''_{2j-1} = f_{m_{2j-1}}$, то
$z'_j=f_{m_{2j-1}}v_j$.

Если же $f''_{2j-1}= f_{m_{2j-1}}^{(k)}$, $k = 1,2$, а подслово
$z'_j$ пересекается с подсловом $z'_{j+1}$, то $f''_{2j}\subset z'_i$. Так
 как слова $x_{(2j)''}$ и $x_{(2j-1)''}$ -- нециклические, то
$|x_{(2j)''}| = |x_{(2j-1)''}|$. Но тогда длина периода в $z'_j$
не меньше $4n$, что противоречит замечанию \ref{c:pr}.

Тем самым доказана следующая лемма:
\begin{lemma}\label{lThick}
В слове $W$ с высотой не более $(\Psi(n,4n,l)+ 5s')$ найдётся не менее $s'$ непересекающихся периодических подслов, в которых период повторится не менее $2n$ раз. Кроме того, между любыми двумя элементами данного множества периодических подслов найдётся подслово длины периода более левого из выбранных элементов.
\end{lemma}

\subsection{Завершение доказательств теоремы \ref{t1:log_2} и основной теоремы~\ref{c:main2}}
Подставляя в лемму \ref{lThick} вместо числа $s'$ значение $s$ из доказательства теоремы \ref{ThThick} получаем, что
 высота $W$ не больше, чем
$$\Psi(n,4n,l)+ 5s < E_1 l\cdot n^{E_2+12\log_3 n} ,$$ где
$E_1 = 4^{21\log_3 4 + 17}, E_2 = 30\log_3 4 + 10.$

Тем самым мы получили {\bf утверждение
основной теоремы \ref{c:main2}}.

{\bf Доказательство теоремы \ref{t1:log_2}} завершается также, только в пункте \ref{ThThick:end} вместо последовательности
$$a_0 = 3^{\ulcorner \log_3 p_{n,d}\urcorner}, a_1 = 3^{\ulcorner \log_3 p_{n,d}\urcorner-1},\ldots,a_{{\ulcorner \log_3 p_{n,d}\urcorner}} =1$$
рассматривается последовательность
$$a_0 = 2^{\ulcorner \log_2 p_{n,d}\urcorner}, a_1 = 2^{\ulcorner \log_2 p_{n,d}\urcorner-1},\ldots,a_{{\ulcorner \log_2 p_{n,d}\urcorner}} =1,$$
а значение $\Psi(n,4n,l)$ берётся из теоремы \ref{t2:log_2}.

\section{Комментарии}

Представленная вниманию читателя техника, по всей видимости, позволяет улучшить
полученную в данную работе оценку, но при этом она останется только
субэкспоненциальной. Для получения полиномиальной оценки, если она
существует, требуются новые идеи и методы.

В начале представленного решения при использовании теоремы Ширшова
подслова большого слова используются прежде всего в качестве
множества независимых элементов, а не набор тесно связанных друг с
другом слов. Далее используется то, что буквы внутри
подслов раскрашены. При учёте раскраски только первых букв
подслов получается экспоненциальная оценка. При рассмотрении
раскраски всех букв подслов опять получается экспонента. Данный
факт имеет место из-за построения иерархической системы подслов. Не исключено,
что подробное рассмотрение приведенной связи подслов вкупе с
изложенным выше решением позволит улучшить полученную оценку
вплоть до полиномиальной.

Интересно также получить оценки на высоту алгебры над множеством
слов степени не выше сложности алгебры (в англоязычной литературе $\PI$-degree). В работе
\cite{BBL97} получены экспоненциальные оценки, а для
слов, не являющихся линейной комбинацией лексикографически меньших,
в работе \cite{Bel07} получены надэкспоненциальные оценки.

Глубокие идеи оригинальных работ А.~И.~Ширшова
\cite{Shirshov1,Shirshov2}, восходящие к технике элиминации в
алгебрах Ли, могут оказаться чрезвычайно полезными, в том числе и для улучшения оценок, несмотря на то,
что оценки на высоту, полученные в этих работах, являются только
примитивно рекурсивными.

\end{document}